\documentclass[fleqn]{mat01}
\usepackage{times,mathtimy,amssymb,latexsym}
\begin{document}

\setcounter{page}{257} \firstpage{257}

\newtheorem{theore}{\bf Theorem}
\newtheorem{defini}{\rm DEFINITION}

\newtheorem{theoree}{\bf Theorem}
\renewcommand\thetheoree{\arabic{section}.\arabic{theoree}}
\newtheorem{theor}[theoree]{\bf Theorem}
\newtheorem{pot}[theoree]{\it Proof of Theorem}

\def\rema{\trivlist \item[\hskip \labelsep{\it Remark.}]}
\def\apoc{\trivlist \item[\hskip \labelsep{\it Another proof of {\rm (}i{\rm )} {\rm (}Combinatorial{\rm )}.}]}

\renewcommand{\theequation}{\arabic{section}.\arabic{equation}}

\title{$\pmb{n}$-Colour self-inverse compositions}

\markboth{Geetika Narang and A~K~Agarwal}{$n$-Colour self-inverse
compositions}

\author{GEETIKA NARANG and A~K~AGARWAL}

\address{Centre for Advanced Studies in Mathematics, Panjab University,
Chandigarh~160~014, India\\
\noindent E-mail: geetika2narang@yahoo.com; aka@pu.ac.in}

\volume{116}

\mon{August}

\parts{3}

\pubyear{2006}

\Date{MS received 18 November 2005; revised 12 April 2006}

\begin{abstract}
MacMahon's definition of self-inverse composition is extended to
$n$-colour self-inverse composition. This introduces four new
sequences which satisfy the same recurrence relation with
different initial conditions like the famous Fibonacci and Lucas
sequences. For these new sequences explicit formulas, recurrence
relations, generating functions and a summation formula are
obtained. Two new binomial identities with combinatorial meaning
are also given.
\end{abstract}


\maketitle

\section{Introduction}

In the classical theory of partitions, compositions were first
defined by MacMahon \cite{6} as ordered partitions. For example,
there are 5 partitions and 8 compositions of 4. The partitions are
$4, 31, 22, 21^{2}, 1^{4}$ and the compositions are $4, 31, 13,
22, 21^{2}, 121,$\break $1^{2}2, 1^{4}$.

Agarwal and Andrews \cite{4} defined an $n$-colour partition as a
partition in which a part of size $n$ can come in $n$ different
colours. We use a subscript to index the colour and write
$n_{1},n_{2},\dots,n_{n}$ for the parts of size $n$. Analogous to
MacMahon's ordinary compositions Agarwal \cite{2} defined an
$n$-colour composition as an $n$-colour ordered partition. Thus,
for example, there are 8 $n$-colour compositions of 3, viz.,
\begin{align*}
&3_{1},~3_{2},~3_{3},\\[.2pc]
&2_{1}1_{1},~2_{2}1_{1},~1_{1}2_{1},~1_{1}2_{2},\\[.2pc]
&1_{1}1_{1}1_{1}.
\end{align*}
More properties of $n$-colour compositions were found in \cite{3}.

\begin{defini}{\rm \cite{6}}$\left.\right.$\vspace{.5pc}

\noindent {\rm A composition is said to be self-inverse when the
parts of the composition read from left to right are identical
with when read from right to left.}
\end{defini}

Analogous to the above definition of classical self-inverse
compositions we define an $n$-colour self-inverse composition as
follows:\pagebreak

\begin{defini}$\left.\right.$\vspace{.5pc}

\noindent {\rm An $n$-colour composition whose parts read from
left to right are identical with when read from right to left is
called an $n$-colour self-inverse composition. Thus, for example,
there are 11 $n$-colour self-inverse compositions of 5 viz.,
\begin{align*}
&5_{1}, 5_{2}, 5_{3}, 5_{4}, 5_{5},\\[.2pc]
&1_{1}3_{1}1_{1}, 1_{1}3_{2}1_{1}, 1_{1}3_{3}1_{1},\\[.2pc]
&2_{1}1_{1}2_{1}, 2_{2}1_{1}2_{2},\\[.2pc]
&1_{1}1_{1}1_{1}1_{1}1_{1}.
\end{align*}}
\end{defini}

\begin{defini}$\left.\right.$\vspace{.5pc}

\noindent {\rm The Lucas sequence $\{L_{n}\}_{n=0}^{\infty}$ is
defined as $L_{0}=2,~L_{1}=1$ and $L_{n}=L_{n-1}+L_{n-2}$ when
$n\geq2$.}
\end{defini}

\begin{defini}$\left.\right.$\vspace{.5pc}

\noindent {\rm A point whose co-ordinates are integers is called a
lattice point (unless otherwise stated, we take these integers to
be non-negative).}
\end{defini}

\begin{defini}$\left.\right.$\vspace{.5pc}

\noindent {\rm By a lattice path (or simply a path), we mean a
minimal path via lattice points taking unit horizontal and unit
vertical steps.}
\end{defini}

\begin{defini}$\left.\right.$\vspace{.5pc}

\noindent {\rm The Fibonacci sequence $\{F_{n}\}_{n=0}^{\infty}$
is defined as $F_{0}=0, F_{1}=1$ and $F_{n}=F_{n-1}+F_{n-2}$ when
$n\geq2$.}
\end{defini}

Agarwal \cite{2} proved the following theorem.

\begin{theore}[\!] Let $C(m,q)$ and $C(q)$ denote the
enumerative generating functions for $C(m,\nu)$ and $C(\nu)${\rm
,} respectively{\rm ,} where $C(m,\nu)$ is the number of
$n$-colour compositions of $\nu$ into $m$ parts and $C(\nu)$ is
the number of $n$-colour compositions of $\nu$. Then
\begin{align}
&C(m,q) = \frac{q^{m}}{(1-q)^{2m}},\\[.2pc]
&C(q) = \frac{q}{1-3q+q^{2}},\\[.3pc]
&C(m,\nu) = \left(\begin{array}{c}
\nu+m-1\\[.2pc]
2m-1
\end{array}\right),\\[.2pc]
&C(\nu) = F_{2\nu}.
\end{align}
\end{theore}

\section{Explicit formulas}

\setcounter{equation}{0}

In this section we shall prove the following explicit formulas for
$n$-colour self-inverse compositions.

\begin{theor}[\!]
Let $A(m,\nu)$ be the number of $n$-colour self-inverse
compositions of $\nu$ into $m$ parts. Then
\begin{align}
\hskip -4pc ({\rm i})\hskip 2.6pc &A(2m-1,2\nu-1)=
\sum_{l=1}^{\nu-1} (2l-1) \left(\begin{array}{c}
\nu+m-l-2\\[.2pc]
2m-3
\end{array}\right)\quad {\rm for} \ \ m>1,\\[.3pc]
\hskip -4pc ({\rm ii})\hskip 2.6pc &A(2m,2\nu) =
\left(\begin{array}{c}
\nu+m-1\\[.2pc]
2m-1
\end{array}\right),\\[.3pc]
\hskip -4pc ({\rm iii})\hskip 2.6pc &A(2m-1,2\nu) =
\sum_{l=1}^{\nu-1} 2l\left(\begin{array}{c}
\nu+m-l-2\\[.2pc]
2m-3
\end{array}\right)\quad {\rm for} \ \ m>1.
\end{align}
\end{theor}

\setcounter{theoree}{0}
\begin{pot}{\rm
We shall prove (i) and the other parts can be proved similarly.
(i)~Obviously, an odd number can have self-inverse $n$-colour
compositions only when the number of parts is odd. A~self-inverse
$n$-colour composition of an odd number $2\nu-1$ into $2m-1~(m>1)$
parts can be read as a central part, say, $2l-1$, and two
identical $n$-colour compositions of $\nu-l$ into $m-1$ parts on
each side of the central part.

The number of $n$-colour compositions of $\nu-l$ into $m-1$ parts
is $\left(\begin{smallmatrix} \nu+m-l-2\\[.1pc]
2m-3 \end{smallmatrix}\right)$ by (1.3). Now, the central part can
appear in $2l-1$ ways and takes values from 1 to $2\nu-3$. So, $l$
varies from 1 to $\nu-1$.

Therefore, for $m>1$, the number of $n$-colour self-inverse
compositions of $2\nu-1$ into $2m-1$ parts is
\begin{equation*}
A(2m-1,2\nu-1) = \sum_{l=1}^{\nu-1} (2l-1)\left(\begin{array}{c}
\nu+m-l-2\\[.2pc] 2m-3 \end{array}\right).
\end{equation*}
For $m=1$, there will be $2\nu-1$ $n$-colour self-inverse
compositions. Therefore, the total number of $n$-colour
self-inverse compositions of $2\nu-1~$ is
\begin{equation*}
(2\nu-1) + \sum_{m=2}^{\nu} \sum_{l=1}^{\nu-1}
(2l-1)\left(\begin{array}{c} \nu+m-l-2\\[.2pc]
2m-3\end{array}\right).
\end{equation*}}
\end{pot}
\vspace{6pt}
\section{Recurrence formulas}

\setcounter{equation}{0}

\setcounter{theoree}{0}

In this section, we shall prove the following: If
$a_{n},~b_{n},~c_{n}$ and $d_{n}$ denote the number of $n$-colour
self-inverse compositions of $2n+1$, $2n$ into an even number of
parts, $2n$ into an odd number of parts and $2n$ (that is,
$d_{n}=b_{n}+c_{n}$) respectively, then the following theorem
holds.

\begin{theor}[\!] The sequences $\{a_{n}\}, \{b_{n}\},\{c_{n}\},
\{d_{n}\}$ satisfy the following recurrences.

\begin{enumerate}
\renewcommand\labelenumi{\rm (\roman{enumi})}
\leftskip .45pc
\item $a_{0}=1, a_{1}=4$ and $a_{n}=3a_{n-1}-a_{n-2}~$ for $n\geq2,$

\item $b_{1}=1, b_{2}=3$ and $b_{n}=3b_{n-1}-b_{n-2}~$ for $n>2,$

\item $c_{1}=2, c_{2}=6$ and $c_{n}=3c_{n-1}-c_{n-2}~$ for $n>2,$

\item $d_{1}=3, d_{2}=9$ and $d_{n}=3d_{n-1}-d_{n-2}~$ for $n>2$ where $d_{n}=b_{n}+c_{n}.$
\end{enumerate}
\end{theor}

\begin{rema}
Theorem~2.1 implies that
\begin{align*}
b_{\nu} = \sum_{m=1}^{\nu} \left(\begin{array}{c}
\nu+m-1\\[.2pc]
2m-1
\end{array}\!\right)\quad\!\hbox{and}\quad\!
c_{\nu} = 2\nu+ \sum_{m=2}^{\nu} \sum_{l=1}^{\nu-1} 2l
\left(\begin{array}{c}
\nu+m-l-2\\[.2pc]
2m-3
\end{array}\right).
\end{align*}
\end{rema}

\setcounter{theoree}{0}
\begin{pot}{\rm (i) We have
\begin{align*}
\hskip -4pc a_{n} &= (2n+1) + \sum_{m=2}^{n+1} \sum_{l=1}^{n}
(2l-1)\left(\begin{array}{c}
n+m-l-1\\[.1pc] 2m-3
\end{array}\right)\\[.3pc]
\hskip -4pc &= (2n+1)+ \sum_{m=2}^{n}
\sum_{l=1}^{n-1}(2l-1)\left(\begin{array}{c}
n+m-l-1\\[.1pc] 2m-3
\end{array}\right)\\[.3pc]
\hskip -4pc &\quad\, + \sum_{m=2}^{n}(2l-1)\left(\begin{array}{c}
n+m-l-1\\[.1pc] 2m-3
\end{array}\right)
\bigg|_{l=n}\\[.3pc]
\hskip -4pc &\quad\, + \sum_{l=1}^{n} (2l-1)
\left(\begin{array}{c}
n+m-l-1\\[.1pc] 2m-3
\end{array}\right)\bigg|_{m=n+1}\\[.3pc]
\hskip -4pc &=(2n+1) + \sum_{m=2}^{n} \sum_{l=1}^{n-1} (2l-1)
\left(\begin{array}{c}
n+m-l-1\\[.1pc] 2m-3
\end{array}\right) + (2n-1)+1\\[.3pc]
\hskip -4pc &= 4n+1+\sum_{m=2}^{n} \sum_{l=1}^{n-1} (2l-1)
\left(\begin{array}{c}
n+m-l-2\\[.1pc] 2m-3
\end{array}\right)\\[.3pc]
\hskip -4pc &\quad\, + \sum_{m=2}^{n} \sum_{l=1}^{n-1}
(2l-1)\left(\begin{array}{c}
n+m-l-2\\[.1pc] 2m-4
\end{array}\right)\\[.3pc]
\hskip -4pc &= 4n+1+\{a_{n-1}-(2n-1)\} + \sum_{m=1}^{n-1}
\sum_{l=1}^{n-1} (2l-1) \left(\begin{array}{c}
n+m-l-1\\[.1pc]
2m-2
\end{array}\right)\\[.3pc]
\hskip -4pc &=2n+2+a_{n-1}+ \sum_{m=2}^{n-1} \sum_{l=1}^{n-1}
(2l-1) \left(\begin{array}{c}
n+m-l-1\\[.1pc] 2m-2
\end{array}\right) + \sum_{l=1}^{n-1}(2l-1)\\[.3pc]
\hskip -4pc &=2n+2+a_{n-1}+ \sum_{m=2}^{n-1}
\sum_{l=1}^{n-1}(2l-1) \left(\begin{array}{c}
n+m-l-2\\[.1pc] 2m-3 \end{array}\right)\\[.3pc]
\hskip -4pc &\quad\, + \sum_{m=2}^{n-1} \sum_{l=1}^{n-1} (2l-1)
\left(\begin{array}{c} n+m-l-2\\[.1pc] 2m-2 \end{array}\right) + \sum_{l=1}^{n-1} (2l-1)\\[.3pc]
\hskip -4pc &\quad\, \left(\hbox{\rm by the
binomial identity}\ \left(\begin{array}{c} n\\[.1pc] m \end{array} \!\right)
= \left(\begin{array}{c} n-1\\[.1pc] m-1 \end{array} \!\right) +
\left(\begin{array}{c} n-1\\[.1pc] m \end{array}\!\right)\right)\\
\hskip -4pc
\end{align*}
\pagebreak
\begin{align*}
\hskip -4pc &= 2n+2+a_{n-1}+\{a_{n-1}-(2n-1)-1\}\\[.3pc]
\hskip -4pc &\quad\, + \sum_{m=2}^{n-1} \sum_{l=1}^{n-2} (2l-1)
\left(\begin{array}{c} n+m-l-2\\ 2m-2 \end{array}\right) \quad
+(n-1)^{2}\\[.3pc]
\hskip -4pc \phantom{a_{n}} &= 2+2a_{n-1}+(n-1)^{2}\\[.3pc]
\hskip -4pc &\quad\, + \sum_{m=2}^{n-1} \sum_{l=1}^{n-2} (2l-1)
\left\{\left(\begin{array}{c} n+m-l-1\\[.1pc]
2m-2
\end{array}\right) - \left(\begin{array}{c} n+m-l-2\\[.1pc]
2m-3
\end{array}\right)\right\}\\[.3pc]
\hskip -4pc &= 2+2a_{n-1}+(n-1)^{2}+\sum_{m=2}^{n-1}
\sum_{l=1}^{n-2} (2l-1) \left(\begin{array}{c} n+m-l-1\\[.1pc] 2m-2
\end{array}\right)\\[.3pc]
\hskip -4pc &\quad\, - \sum_{m=2}^{n-1} \sum_{l=1}^{n-2} (2l-1)
\left(\begin{array}{c} n+m-l-3\\[.1pc]
2m-3
\end{array}\right) - \sum_{m=2}^{n-1}
\sum_{l=1}^{n-2} (2l-1) \left(\begin{array}{c} n+m-l-3\\[.1pc] 2m-4
\end{array}\right)\\[.3pc]
\hskip -4pc &= 2+2a_{n-1}+(n-1)^{2}+ \left\{ \sum_{m=2}^{n-1}
\sum_{l=1}^{n-2} (2l-1) \left(\begin{array}{c} n+m-l-2\\[.1pc] 2m-2
\end{array}\right)\right.\\[.3pc]
\hskip -4pc &\quad\, + \sum_{m=2}^{n-1} \sum_{l=1}^{n-2} (2l-1)
\left. \left(\begin{array}{c} n+m-l-2\\[.1pc] 2m-3 \end{array}
\right)\right\} - \{a_{n-2} - (2n-3)\}\\[.3pc]
\hskip -4pc &\quad\, - \sum_{m=2}^{n-1} \sum_{l=1}^{n-2} (2l-1)
\left(\begin{array}{c} n+m-l-3\\[.1pc] 2m-4
\end{array}\right)\\[.3pc]
\hskip -4pc &= 2a_{n-1} - a_{n-2} + (n-1)^{2} + (2n-1) +
\sum_{m=2}^{n-1} \sum_{l=1}^{n-2} (2l-1)\left(\begin{array}{c}
n+m-l-2\\[.1pc]
2m-3
\end{array}\right)\\[.3pc]
\hskip -4pc &\quad\, + \sum_{m=2}^{n-1} \sum_{l=1}^{n-2}
(2l-1)\left(\begin{array}{c}
n+m-l-2\\[.1pc] 2m-2 \end{array}\right) - \sum_{m=2}^{n-1}
\sum_{l=1}^{n-2} (2l-1) \left(\begin{array}{c}
n+m-l-3\\[.1pc] 2m-4 \end{array}\right)\\[.3pc]
\hskip -4pc &= 2a_{n-1}-a_{n-2}+(n-1)^{2}+(2n-1)\\[.3pc]
\hskip -4pc &\quad\, + \left\{ \sum_{m=2}^{n} \sum_{l=1}^{n-1}
(2l-1) \left(\begin{array}{c} n+m-l-2\\[.1pc]
2m-3
\end{array}\right) - 1 - (2n-3) \right\}\\[.4pc]
\hskip -4pc &\quad\, + \left\{ \sum_{m=2}^{n-2} \sum_{l=1}^{n-2}
(2l-1) \left(\begin{array}{c} n+m-l-2\\[.1pc] 2m-2 \end{array}\right) +
\sum_{l=1}^{n-2} (2l-1) \left(\begin{array}{c} 2n-l-3\\[.1pc] 2n-4
\end{array}\right)\right\}\\[.4pc]
\hskip -4pc &\quad\, -\left\{ \sum_{m=2}^{n-2} \sum_{l=1}^{n-2}
(2l-1) \left(\begin{array}{c} n+m-l-2\\[.1pc] 2m-2 \end{array} \right) +
\sum_{l=1}^{n-2} (2l-1) \left(\begin{array}{c} n-l-1\\[.1pc] 0
\end{array}\right) \right\}\\[.35pc]
\hskip -4pc &=2a_{n-1}-a_{n-2}+(n-1)^{2}+(2n-1)+
\{a_{n-1}-(2n-1)-1-(2n-3)\}\\[.3pc]
\hskip -4pc &\quad\, +1-(n-2)^{2}%\\[.3pc]
\end{align*}
\begin{align*}
\hskip -4pc \phantom{a_{n}} &=3a_{n-1}-a_{n-2}+(n-1)^{2}-(2n-3)-(n-2)^{2}\\[.3pc]
\hskip -4pc &=3a_{n-1}-a_{n-2}.
\end{align*}
The other parts can be proved on the same lines.}
\end{pot}

\begin{apoc}
To prove that $a_{n}+a_{n-2}=3a_{n-1},$ we split the
\hbox{$n$-colour} compositions enumerated by $a_{n}+a_{n-2}$ into
three classes:\vspace{-.3pc}
\begin{enumerate}
\renewcommand\labelenumi{(\roman{enumi})}
\leftskip .4pc
\item enumerated by $a_{n}$ and having $1_{1}$ on both extremes,

\item enumerated by $a_{n}$ and having $h_{t}$ on both extremes,
$h>1$ and $1\leq t \leq h-1$ and $n$-colour self-inverse
compositions of $2n+1$ of the form $(2n+1)_{t}, ~1\leq t \leq
2n-1.$

\item enumerated by $a_{n}$ and having $h_{h}$ on both extremes,
$h>1$, $(2n+1)_{2n}, (2n+1)_{2n+1}$ and those enumerated by
$a_{n-2}$.\vspace{-.5pc}
\end{enumerate}
\end{apoc}
We transform the $n$-colour compositions in class~(i)~by deleting
$1_{1}$ on both extremes. This produces an $n$-colour composition
enumerated by $a_{n-1}$. Conversely, given any $n$-colour
composition enumerated by $a_{n-1}$ we add $1_{1}$ on both
extremes to produce the elements of the first class. In this way
we establish that there are exactly $a_{n-1}$ elements in the
first class. Next, we transform the $n$-colour compositions in
class~(ii)~by subtracting~1 from $h$, that is, replacing $h_{t}$
by $(h-1)_{t}$ and subtracting~2 from $(2n+1)$ of $(2n+1)_{t},
~1\leq t \leq 2n-1$. This transformation also establishes the fact
that there are exactly $a_{n-1}$ elements in class~(ii). Finally,
we transform the elements in class (iii) as follows: Subtract
$1_{1}$ from $h_{h}$, $h<2n+1$, that is, replace $h_{h}$ by
$(h-1)_{h-1}$. We will get those $n$-colour self-inverse
compositions of $2(n-1)+1$ whose extremes are $h_{h}$, except
self-inverse compositions in one part only. Also replace
$(2n+1)_{2n}$ by $(2n-1)_{2n-2}$ and $(2n+1)_{2n+1}$ by
$(2n-1)_{2n-1}$, and to get the remaining $n$-colour self-inverse
compositions from $a_{n-2}$ we add 1 to both extremes, that is,
replace $h_{t}$ by $(h+1)_{t}$. For $n$-colour self-inverse
compositions into one part we add 2, that is, replace $(2n-3)_{t}$
by $(2n-1)_{t}$, $1 \leq t \leq 2n-3$. We see that the number of
$n$-colour compositions in class (iii) is also equal to $a_{n-1}$.

So each class has $a_{n-1}$ elements. Hence,
$a_{n}+a_{n-2}=3a_{n-1}$.

\section{Generating functions}

\setcounter{equation}{0}

\setcounter{theoree}{0}

In this section, we will prove the following theorem.

\begin{theor}[\!]
\begin{align*}
\hskip -4pc {\rm
(i)} \hskip 2.6pc
&\sum_{n=0}^{\infty}a_{n}q^{n}=\frac{1+q}{q^{2}-3q+1},\\[.3pc]
\hskip -4pc {\rm
(ii)} \hskip 2.6pc &\sum_{n=1}^{\infty}b_{n}q^{n}=\frac{q}{q^{2}-3q+1},\\[.3pc]
\hskip -4pc {\rm (iii)} \hskip 2.6pc
&\sum_{n=1}^{\infty}c_{n}q^{n}=\frac{2q}{q^{2}-3q+1},~{\rm
and}\\[.3pc]
\hskip -4pc {\rm (iv)} \hskip 2.6pc
&\sum_{n=1}^{\infty}d_{n}q^{n}=\frac{3q}{q^{2}-3q+1}.
\end{align*}
\end{theor}

\setcounter{theoree}{0}
\begin{pot}{\rm
\begin{align*}
\hskip -4pc {\rm (i)} \hskip 8pc \sum_{n=0}^{\infty} a_{n}q^{n}
&= a_{0}+ a_{1}q +\sum_{n=2}^{\infty}a_{n}q^{n}\\[.2pc]
&=1+4q+\sum_{n=2}^{\infty}(3a_{n-1}-a_{n-2})q^{n}\\[.2pc]
&=1+4q+3\sum_{n=1}^{\infty}a_{n}q^{n+1}-
\sum_{n=0}^{\infty}a_{n}q^{n+2}\\[.2pc]
&=1+4q+3\sum_{n=0}^{\infty}a_{n}q^{n+1}
-3q-\sum_{n=0}^{\infty}a_{n}q^{n+2}\\[.2pc]
\hskip -4pc (q^{2}-3q+1)\sum_{n=0}^{\infty}a_{n}q^{n} &=1+q\\[.2pc]
\hskip -4pc \sum_{n=0}^{\infty}a_{n}q^{n}&=\frac{1+q}{q^{2}-3q+1}.
\end{align*}
This completes the proof of (i). The proofs of (ii) and (iii) are
similar and hence are omitted. By adding (ii) and (iii) we get
(iv).}
\end{pot}

\section{A summation formula}

\setcounter{equation}{0}

\setcounter{theoree}{0}

\begin{theor}[\!]
For $n \geq 0${\rm ,}
\begin{equation*}
\sum_{k=0}^{n}(-1)^{n+k} \left(\begin{array}{c} 2n+1\\ n-k
\end{array} \right) a_{k}=1,
\end{equation*}
where
\begin{equation*}
a_{k}=(2k+1) + \sum_{m=2}^{k+1} \sum_{l=1}^{k} (2l-1)
\left(\begin{array}{c} k+m-l-1\\ 2m-3 \end{array} \right).
\end{equation*}
\end{theor}

\setcounter{theoree}{0}
\begin{pot}{\rm
We will prove the result by induction on $n$. The result is
obviously true for $n= 0,1$. We assume that the result is true for
$n$, that is,
\begin{equation*}
\sum_{k=0}^{n}(-1)^{n+k} \left(\begin{array}{c} 2n+1\\[.1pc] n-k
\end{array} \right) a_{k}=1.
\end{equation*}
Now
\begin{align*}
&\sum_{k=0}^{n+1} (-1)^{n+1+k} \left(\begin{array}{c} 2n+3\\[.1pc]
n+1-k
\end{array}\right)
a_{k}\\[.2pc]
&\quad\, = \sum_{k=0}^{n+1}(-1)^{n+1+k} \left\{
\left(\begin{array}{c}
2n+2\\[.1pc] n+1-k\end{array} \right)+ \left(\begin{array}{c} 2n+2\\[.1pc]
n-k\end{array} \right) \right\}
a_{k}%\\[.2pc]
\end{align*}
\begin{align*}
&\quad\, = \sum_{k=0}^{n+1}(-1)^{n+1+k}\left\{
\left(\begin{array}{c}
2n+1\\[.1pc] n+1-k\end{array}\right)\right. + 2
\left(\begin{array}{c} 2n+1\\[.1pc] n-k\end{array}\!\right) +
\left.\left(\begin{array}{c} 2n+1\\[.1pc]
n-k-1 \end{array}\!\right)\right\}a_{k}\\[.2pc]
&\qquad\ \left({\rm by}\ \left(\begin{array}{c} n\\[.1pc] k\end{array}\right)=
\left(\begin{array}{c} n-1\\[.1pc] k\end{array}\right) +
\left(\begin{array}{c} n-1\\[.1pc] k-1 \end{array}\right)\right)\\[.2pc]
&\quad\, = \sum_{k=0}^{n+1}(-1)^{n+1+k} \left(\begin{array}{c} 2n+1\\[.1pc]
n+1-k\end{array}\right) a_{k}-2 +
\sum_{k=1}^{n}(-1)^{n+k} \left(\begin{array}{c} 2n+1\\[.1pc] n-k\end{array}\right) a_{k-1}\\[.2pc]
&\quad\, = \sum_{k=0}^{n+1}(-1)^{n+1+k} \left(\begin{array}{c} 2n+1\\[.1pc]
n+1-k \end{array} \right)a_{k}-2+3-3(-1)^{n} \left(\begin{array}{c} 2n+1\\[.1pc] n
\end{array}\right)\\[.2pc]
&\qquad\ -\sum_{k=2}^{n+1}(-1)^{n+1+k}
\left(\begin{array}{c} 2n+1\\[.1pc]
n+1-k
\end{array} \right) a_{k}\\[.2pc]
&\quad\, = \sum_{k=0}^{1}(-1)^{n+1+k} \left(\begin{array}{c} 2n+1\\[.1pc]
n+1-k \end{array}\right) a_{k}+1-3(-1)^{n} \left(\begin{array}{c}
2n+1\\[.1pc] n\end{array}\right)\\[.2pc]
&\quad\, =1.
\end{align*}
Therefore the result is true for $n+1$.}
\end{pot}

Hence the result is true for all $n\in \{0,1,2,\dots\}$.

\section{New binomial identities with combinatorial interpretations}

\setcounter{equation}{0}

\setcounter{theoree}{0}

\begin{theor}[\!]
For $\nu\geq0${\rm ,}
\begin{align*}
\hskip -4pc \hbox{\rm (i)} \hskip 2.95pc &(2\nu+1) +
\sum_{m=2}^{\nu+1} \sum_{l=1}^{\nu}(2l-1) \left(\begin{array}{c}
\nu + m-l-1\\[.1pc] 2m-3
\end{array}\right)\\[.2pc]
\hskip -4pc &\quad\, = \sum_{m=0}^{\nu} \frac{2\nu+1}{2\nu+1-m}
\left(\begin{array}{c} 2\nu+1-m\\[.1pc] m \end{array}\right),\\[.2pc]
\hskip -4pc \hbox{\rm (ii)} \hskip 2.95pc &\nu + \sum_{m=2}^{\nu}
\sum_{l=1}^{\nu-1} l\left(\begin{array}{c} \nu + m-l-2\\[.1pc] 2m-3
\end{array} \right) = \sum_{m=1}^{\nu} \left(\begin{array}{c} \nu
+ m-1\\[.1pc] 2m-1
\end{array} \right).
\end{align*}
\end{theor}

\setcounter{theoree}{0}
\begin{pot}{\rm
$\left.\right.$\vspace{.5pc}

\noindent (i) Let
\begin{equation*}
A_{\nu}= (2\nu+1)+ \sum_{m=2}^{\nu+1} \sum_{l=1}^{\nu} (2l-1)
\left(\begin{array}{c} \nu + m-l-1\\ 2m-3
\end{array}\right)
\end{equation*}
and
\begin{equation*}
B_{\nu} = \sum_{m=0}^{\nu} \frac{2\nu+1}{2\nu+1-m}
\left(\begin{array}{c} 2\nu+1-m\\ m \end{array}\right).
\end{equation*}
Clearly, $A_{0}=B_{0}=1, ~A_{1}=B_{1}=4,$ and we have already
proved $A_{\nu}=3A_{\nu-1}-A_{\nu-2}$ for $\nu\geq2$.\pagebreak

Now, for $\nu\geq2$, we have
\begin{align*}
\hskip -4pc B_{\nu} &= \sum_{m=0}^{\nu} \frac{2\nu+1}{2\nu+1-m}
\left(\begin{array}{c} 2\nu + 1 - m\\[.1pc] m \end{array}\right)\\[.2pc]
\hskip -4pc &= \sum_{m=0}^{\nu} \left\{\left(\begin{array}{c} 2\nu + 1 - m\\[.1pc]
m \end{array}\right) + \left(\begin{array}{c} 2\nu-m\\[.1pc] m-1
\end{array}\right)\right\}\\[.2pc]
\hskip -4pc &\quad\, \left(\hbox{\rm by the binomial identity}\
\frac{n}{n-m}\left(\begin{array}{c} n-m\\[.1pc] m
\end{array}\right) = \left(\begin{array}{c} n - m\\[.1pc] m
\end{array}\right) + \left(\begin{array}{c} n-m-1\\[.1pc] m-1
\end{array} \right) \right)\\[.2pc]
\hskip -4pc &= \sum_{m=0}^{\nu} \left\{ \left(\begin{array}{c} 2\nu-m\\[.1pc]
m \end{array} \right)+ \left(\begin{array}{c} 2\nu-m-1\\[.1pc] m-1
\end{array}\right)\right\} + B_{\nu-1}\\[.2pc]
\hskip -4pc &\quad\, \left(\hbox{\rm by the binomial identity}\
\left(\begin{array}{c} n\\[.1pc] m\end{array}\right) \right.\\[.2pc]
\hskip -4pc &\qquad\, \left. =
\left(\begin{array}{c} n-1\\[.1pc] m
\end{array} \right) + \left(\begin{array}{c} n-1\\[.1pc] m-1
\end{array}\right) \hbox{\rm and the definition of}\ B_{\nu}\right)\\[.2pc]
\hskip -4pc &= B_{\nu-1}+B_{\nu-1}+ \sum_{m=0}^{\nu} \left\{
\left(\begin{array}{c} 2\nu-1-m\\[.1pc] m-1 \end{array}\right) +
\left(\begin{array}{c} 2\nu-2-m\\[.1pc] m-2
\end{array}\right)\right\}\\[.2pc]
\hskip -4pc &= 2B_{\nu-1} + \sum_{m=0}^{\nu}
\left\{\left(\begin{array}{c} 2\nu-m\\[.1pc] m-1 \end{array}\right) -
\left(\begin{array}{c} 2\nu-m-1\\[.1pc] m-2 \end{array}\right) +
\left(\begin{array}{c}2\nu-1-m\\[.1pc] m-2\end{array}\right)\right.\\[.2pc]
\hskip -4pc &\quad\, -\left. \left(\begin{array}{c} 2\nu-2-m\\[.1pc] m-3
\end{array}\right)\right\}\\[.2pc]
\hskip -4pc &= 2B_{\nu-1} + \sum_{m=0}^{\nu} \left\{
\left(\begin{array}{c} 2\nu-m\\[.1pc] m-1 \end{array} \right)+
\left(\begin{array}{c} 2\nu-1-m\\[.1pc] m-2 \end{array}\right) \right\}\\[.2pc]
\hskip -4pc &\quad\, - \sum_{m=0}^{\nu} \left\{
\left(\begin{array}{c}
2\nu-m-1\\[.1pc] m-2 \end{array} \right) + \left(\begin{array}{c} 2\nu-2-m\\[.1pc]
m-3\end{array}\right)\right\}\\[.2pc]
\hskip -4pc
&=2B_{\nu-1}+\sum_{m=0}^{\nu-1}\left\{\left(\begin{array}{c}
2\nu-m-1\\[.1pc] m \end{array} \right)+ \left(\begin{array}{c} 2\nu-m-2\\[.1pc]
m-1 \end{array}\right) \right\}\\[.2pc]
\hskip -4pc &\quad\, - \sum_{m=0}^{\nu-2} \left\{
\left(\begin{array}{c} 2\nu-m-3\\[.1pc] m \end{array}\right) +\left(\begin{array}{c}2 \nu-m-4\\[.1pc] m-1
\end{array}\right) \right\}\\[.2pc]
\hskip -4pc &= 2B_{\nu-1}+B_{\nu-1}-B_{\nu-2}\\[.2pc]
\hskip -4pc &= 3B_{\nu-1}-B_{\nu-2}.
\end{align*}
Hence $A_{\nu}=B_{\nu}$, for all $\nu\geq0$ as these satisfy the
same initial conditions and the same recurrence relation.

\noindent (ii) By comparing generating functions for $b_{n}$ and
$c_{n}$ in \S4 and considering the remark below Theorem~3.1, we
get (ii) easily.}
\end{pot}

The following theorem provides a combinatorial interpretation of
Theorem~6.1.

\begin{theor}[\!]$\left.\right.$
\begin{enumerate}
\renewcommand\labelenumi{\rm (\alph{enumi})}
\leftskip .1pc
\item For $\nu\geq0${\rm ,} let $A(\nu)$ denote the number of
self-inverse $n$-colour compositions of $~2\nu+1$. Let $B(\nu)
=\sum_{k=0}^{\nu} d_{\nu, k}${\rm ,} where $d_{\nu, k}$ is the
number of lattice paths from $(0,0)$ to $(2\nu+1-k, k)$ under the
following conditions:
\medskip
\leftskip .3pc
\item The paths do not have two vertical steps in succession which
ensures that they do not cross $y=x+1$.

\item Both the first and last steps can not be
vertical.\vspace{-.5pc}
\medskip

Then $A(\nu) = B(\nu)$ for all $\nu$.

\item The number of $n$-colour self-inverse compositions of $2\nu$
equals three times the number of $n$-colour compositions of
$\nu$.\vspace{-.7pc}
\end{enumerate}
\end{theor}

\setcounter{theoree}{1}
\begin{pot}$\left.\right.${\rm
\begin{enumerate}
\renewcommand\labelenumi{\rm (\alph{enumi})}
\leftskip .1pc
\item In Theorem~3.1 we have proved that
$A(\nu)$ satisfies the following recurrence relation
\begin{align*}
\hskip -1.25pc A(0)=1,\quad A(1)=4 \quad \hbox{\rm and}\quad
A(\nu)=3A(\nu-1)-A(\nu-2),\quad \nu \geq 2.
\end{align*}
Agarwal \cite{1} proved that $L_{2\nu+1}$ also satisfy the same
initial conditions and recurrence relation. This shows that
$A(\nu)= L_{2\nu+1}$. He also showed that $L_{2\nu+1} =
\sum_{k=0}^{\nu} d_{\nu, k}$, and it is known that the right-hand
side of (i) of Theorem~6.1 equals $L_{2\nu+1}$ (see p.~185 of
\cite{5}). This leads to the first part of the theorem.

\item From Theorem~4.1, we have$~ d_{\nu}=3b_{\nu}$. Now
(b)~follows from the fact that $b_{\nu}$ also equals the number of
$n$-colour compositions of $\nu$ in view of (1.2) and
Theorem~4(ii).\vspace{-.7pc}
\end{enumerate}}
\end{pot}

\section*{Acknowledgement}

The authors are thankful to the referee for his suggestions which
led to a better presentation of the paper. One of the authors (GN)
is supported by CSIR Award No.~F.No.~9/135(468)/2k3-EMR-I and the
other author (AKA) is supported by CSIR Research Grant
No.~25(0128)/02/EMR-II.

\end{document}